\documentclass[11pt,twoside]{amsart}
\usepackage{amsfonts}
\usepackage{epsfig,graphics}
\usepackage{amssymb}
\usepackage{amscd}

\newtheorem{theorem}{Theorem}[section]

\newtheorem{proposition}[theorem]{Proposition}
       
\newtheorem{definition}[theorem]{Definition}

\newtheorem{lemma}[theorem]{Lemma}

\newtheorem{remark}[theorem]{Remark}

\newcommand{\Hom}{\text{Hom}}
\newcommand{\End}{\text{End }}

\newcommand{\rk}{\text{rk }}

\newcommand{\im}{\text{im }}

\newcommand{\SL}{\text{SL}}
\newcommand{\GL}{\text{GL}}

\newcommand{\Id}{\text{Id}}
\newcommand{\tr}{\text{tr}}

\newcommand{\Gr}{\text{Gr}}

\newcommand{\Aut}{\text{Aut }}

\newcommand{\BR}{{\bf R}}

\newcommand{\lieg}{{\mathfrak{g}}}

\newcommand{\liep}{{\mathfrak{p}}}

\let\x\times
\let\o\circ
\let\Cal\mathcal
\renewcommand{\Im}{\operatorname{Im}}
\newcommand{\La}{\Lambda}
\newcommand{\al}{\alpha}

\newenvironment{Pf}{\medskip \noindent {\bf Proof: }}
   {$\diamondsuit$ }

\begin{document}
\pagenumbering{arabic}
\title[$C^1$ deformations of almost Grassmannian structures]{$C^1$ deformations of almost Grassmannian structures with strongly essential symmetry}
\author{Andreas \v{C}ap and Karin Melnick}
\date{\today}
\thanks{The authors wish to thank Katharina Neusser for helpful conversations.  Melnick gratefully acknowledges support from NSF grant DMS 1255462 and from the Max-Planck-Institut f\"ur Mathematik in Bonn, where she was a Visiting Scientist during much of the writing of this paper.}

\begin{abstract}
We construct a family of $(2,n)$-almost Grassmannian structures of regularity $C^1$, each admitting a one-parameter group of strongly essential automorphisms, and each not flat on any neighborhood of the higher-order fixed point.  This shows that Theorem 1.3 of \cite{mn.1graded} does not hold assuming only $C^1$ regularity of the structure (see also \cite[Prop 3.5]{cap.me.parabolictrans}).
\end{abstract}

\maketitle

\section{Introduction}

Almost-Grassmannian structures belong to the class of \emph{irreducible parabolic geometries} (also called \emph{almost-Hermitian symmetric} structures), which include projective and conformal structures, among many others.  An \emph{$(m,n)$-almost-Grassmannian structure} on an $mn$-dimensional manifold $M$ comprises a vector bundle isomorphism of $TM$ with $\mathcal{E}^* \otimes \mathcal{F}$, where $\mathcal{E}$ and $\mathcal{F}$ are vector bundles over $M$ of respective ranks $m$ and $n$, together with an isomorphism $\wedge^m \mathcal{E} \cong \wedge^n \mathcal{F}$; the latter corresponds to a volume form compatible with the tensor product.  Denote $\Gr(m,n)$ the real Grassmannian variety of $m$-planes in $\BR^{m+n}$, by $\mathcal{E}$ its tautological $m$-plane bundle, and by $\mathcal{F}$ the rank-$n$ anti-tautological bundle.  An $(m,n)$-almost-Grassmannian structure mimics the isomorphism of $T \Gr(m,n)$ with $\mathcal{E}^* \otimes \mathcal{F}$.

Almost-Grassmannian structures have been studied under the guise of  \emph{Segr\'e structures}.   The \emph{Segr\'e cone} $S(m,n)$ is the variety in $\BR^{mn}$ comprising the rank-one elements under the identification with $\Hom(\BR^m, \BR^n)$.
An $(m,n)$-Segr\'e structure on $M^{mn}$ is a bundle of Segr\'e cones $S_x(m,n) \subset T_x M$.  It is essentially equivalent to an $(m,n)$-almost Grassmannian structure (see \cite{akivis.goldberg.ag.structures}).

In the special case $m=n=2$, when $\dim M =4$, an almost-Grassmannian structure
is equivalent to a conformal spin structure of split signature
$(2,2)$. In fact, $(2,n)$-almost-Grassmannian structures in many respects can be
viewed as higher-dimensional analogs of signature-$(2,2)$ conformal
geometry, which is one of the reasons for the interest in them. Torsion-free $(2,n)$-almost-Grassmannian structures correspond to
(anti-)self-dual conformal structures. 

There is a close relation between
$(2,n)$-almost Grassmannian structures and almost quaternionic structures (see
\cite{salamon.quaternionic}).  They can be viewed as different real forms of the same complex
geometry. 

Finally, $(2,n)$-almost-Grassmannian structures are connected to projective
structures by twistor theory via \emph{path geometries}.  The latter are parabolic
geometries that model systems of second-order ODEs.  More precisely, they correspond
to collections of unparametrized curves $\mathcal{C}$ in a manifold $X$ obtained as
the solutions of such a system; these lift to a foliation $\widetilde{\mathcal{C}}$
of the projectivized tangent bundle ${\bf P}(TX)$.
For one special class of path geometries, $\mathcal{C}$ are the unparametrized
geodesics of an affine connection on $M$---that is, a projective structure
,---while for another class, the \emph{path space} ${\bf P}(TX) /
\widetilde{\mathcal{C}}$ locally inherits a $(2,n)$-almost Grassmannian structure 
.  The intersection of the two classes is the \emph{flat} path geometry, for which
$\mathcal{C}$ comprises projective lines in $X = {\bf RP}^{n+1}$ (see \cite[Secs 2
and 3]{grossman.thesis}).

Irreducible parabolic geometries can admit certain very special automorphisms which fix a point and have trivial derivative at that point, which is then called a \emph{higher-order fixed point}.  Note that a semi-Riemannian metric or an affine connection does not admit nontrivial automorphisms with higher-order fixed points.  These \emph{strongly essential automorphisms} occur in abundance on the homogeneous model spaces for each parabolic geometry.  A structure that is locally equivalent to this model is said to be \emph{flat}.  (See 
Sections \ref{model} and \ref{sec.auts} below).  

Many rigidity results say that a strongly essential flow can occur only in the presence of flatness.  Let $x_0$ be a higher-order fixed point.

\begin{itemize}
\item Nagano and Ochiai \cite{nagano.ochiai.proj} proved for a torsion-free connection that existence of a strongly essential projective flow implies projective flatness of the connection on a neighborhood of $x_0$.  

\item  The second author and Neusser proved the analogous result for almost-c-projective structures and almost-quaternionic structures in \cite[Thms 4.4, 1.2]{mn.1graded}. (See also  \cite[Thm 3.7]{cap.me.parabolictrans} for a precursor result on almost-quaternionic structures.)
\end{itemize}

In conformal Lorentzian geometry, Frances smoothly deformed the Minkowski metric in a neighorhood of a point $x_0$ so that it retains a conformal flow with $x_0$ as higher-order fixed point.   The resulting $C^\infty$ metric is conformally flat inside the light cone of $x_0$, but nonflat outside \cite[Sec 6]{frances.ccvf}.  Then came the following rigidity results:

\begin{itemize}
\item In semi-Riemannian geometry, Frances and the second author proved in \cite{fm.champsconfs} that existence of a strongly essential conformal flow implies conformal flatness on an open set $U$ with $x_0 \in \overline{U}$.
\item In \cite[Thm 1.3]{mn.1graded}, the second author and Neusser proved that a $(2,n)$-almost-Grassmannian structure admitting a strongly essential flow is flat on an open set $U$ with $x_0 \in \overline{U}$ (see also \cite[Prop 3.5]{cap.me.parabolictrans} for a partial result). 
\end{itemize}

 Kruglikov and The exhibited a $C^\omega$ homogeneous path geometry which is not flat and admits a strongly essential flow in \cite[Prop 5.3.2]{kruglikov.the.submax}. Path geometries are not irreducible.  The local path space in their example admits a $(2,n)$-almost-Grassmannian structure. The flow descends, but it is not strongly essential on the quotient.

The proofs of the rigidity theorems cited above, as well as the construction of \cite{kruglikov.the.submax}, make use of the \emph{Cartan geometry} canonically associated to the parabolic geometric structures in question.
This association is only possible with sufficiently high regularity; the minimal order required depends on the structure.  

\subsection{Our examples}

In \cite{cap.deformations}, the first author described the infinitesimal automorphisms and deformations of a parabolic geometry intrinsically in terms of the associated Cartan geometry, using the twisted de-Rham sequence of differential froms with coefficients in the adjoint tractor bundle and the corresponding BGG sequence of invariant differential operators. Motivated by this description of infinitesimal deformations, we explicitly construct a family, locally on $\Gr(2,n)$, that is invariant by a strongly essential flow and integrates to a family of deformed structures, all admitting this flow as automorphisms.  These show that Theorem 1.3 of \cite{mn.1graded} does not hold assuming only $C^1$ regularity.

An almost-Grassmannian structure is said to be $C^k$ if $M$ is at least $C^{k+1}$; $\mathcal{E}$, $\mathcal{F}$, and the isomorphism $\wedge^m \mathcal{E} \cong \wedge^n \mathcal{F}$ are at least $C^{k+1}$; and the isomorphisms $TM \cong \mathcal{E}^* \otimes \mathcal{F}$ are $C^k$.  Such structures are said to be equivalent if they are $C^k$ equivalent (see Section \ref{sec.auts} below).

\begin{theorem}
\label{main.thm}
Let $n\geq 3$ and $x_0 \in \Gr(2,n)$.  There are a dense, open neighborhood $U$ of $x_0$; a strongly essential flow $\{z^t\} < \Aut \Gr(2,n)$ with $x_0$ as higher order fixed point; and an $(n-1)$-parameter family of $C^1$ almost-Grassmannian structures of type $(2,n)$ on $U$, of which each:
\begin{itemize}
\item contains $\{z^t \}$ in its automorphism group;
\item is not locally equivalent to $\Gr(2,n)$ on any open set $V$ with $x_0 \in \overline{V}$;
\end{itemize}
  \end{theorem}
  
  The deformations are given in Section \ref{sec.endo.section}, and the precise claims about them are in Proposition \ref{main.prop}.
  
  \begin{remark}
In fact, none of these deformed structures are locally equivalent to the path space of a path geometry; the harmonic torsion is the full obstruction to this property (see \cite[Props 4.4.3, 4.4.45]{cap.slovak.book.vol1}).
  \end{remark}
  
  
  \section{Background}
  
  \subsection{Almost-Grassmannian structures as first-order $G$-structures}\label{AGr}
For almost-Grassmannian structures of low regularity, as we construct below, the description as Cartan geometries is not
available. Thus we start by reviewing various descriptions of such structures with a
special emphasis on the requirements on regularity.

Let us fix integers $m,n\geq 2$ as above,
with the
case $m=2$, $n>2$ being of primary interest. 
An almost-Grassmannian structure as defined above can be equivalently defined as a (first-order)
G-structure for the Lie group
$$
G_0:=\{(A,B)\in\GL(m,\BR)\times\GL(n,\BR):\det(A)\det(B)=1\}. 
$$
Under the identification $\BR^{mn} \cong \Hom(\BR^m,\BR^n)$, 
the natural representation of $G_0$ on $\BR^{mn}$ is 
$\rho(A,B)\cdot X:=BXA^{-1}$. Observe that the resulting homomorphism $G_0\to \GL(mn,\BR)$
has two element kernel $\{(\Id,\Id),(-\Id,-\Id)\}$ and thus is infinitesimally
injective, so this indeed defines a type of
first-order G-structures on
manifolds of dimension $mn$.  

Such a structure is
given by a principal bundle $p_0:\mathcal{G}_0\to M$ with structure group $G_0$
together with a $\rho$-equivariant bundle morphism to the first order frame bundle $\Cal PM$ of $M$.
The structure is $C^k$ provided $\Cal{G}_0$ is a $C^{k+1}$ principal bundle and the morphism to $\Cal PM$ is $C^k$. 

\begin{proposition}\label{prop:G-str}
On a smooth manifold of dimension $mn$, a $C^k$ first-order $G_0$-structure is equivalent
to a $C^k$ almost-Grassmannian structure of type $(m,n)$. 
\end{proposition}
%

The proof is standard; we give the main points. The bundles $\Cal E^*$ and $\Cal F$ are
associated bundles 
to $\Cal G_0$, while conversely $\Cal G_0$ is obtained as a subbundle of the
fibered product of the frame bundles of $\Cal E^*$ and $\Cal F$. This shows that a $\rho$-equivariant bundle morphism
from $\Cal G_0$ to $\Cal P M$ is equivalent to an isomorphism $TM \stackrel{\sim}{\to} \Cal E^*\otimes
\Cal F$, and the correspondence respects $C^k$ regularity.  

A $\rho$-equivariant bundle morphism $\Cal G_0\to\Cal PM$ can be
equivalently encoded as a one-form $\theta\in\Omega^1(\Cal G_0,\BR^{mn})$ 
which is $G_0$--equivariant and strictly horizontal. Denoting by
$r^g:\Cal G_0\to\Cal G_0$ the principal action of $g\in G_0$, equivariance
means $(r^g)^*\theta= \rho(g)^{-1}\o\theta$.
The second condition says that for each point
$u\in\Cal G_0$, the kernel of $\theta(u):T_u\Cal G_0\to\BR^{mn}$ is the vertical
subspace in $T_u \Cal G_0$.
 In this picture, 
 $C^k$ regularity means $\theta$ is $C^k$, in the sense that for each
$C^k$ vector field $\xi\in\frak X(\Cal G_0)$, the function
$\theta(\xi):\Cal G_0\to\BR^{mn}$ is $C^k$.

\subsection{The homogeneous model---the Grassmann variety}\label{model}

In this section we describe the $(m,n)$-almost-Grassmannian structure on $\Gr(m,n)$.

The group $G_0$ can be realized as the subgroup of $G:=\SL(m+n,\BR)$ respecting the decomposition $\BR^{m+n} \cong \BR^m\oplus\BR^n$.
The Lie algebra $\frak g_0$ is identified with the corresponding block diagonal subalgebra of $\lieg$. Its adjoint action on $\lieg$ preserves a decomposition $\lieg=\lieg_{-1}\oplus\lieg_0\oplus\lieg_1$, where $\lieg_{-1}$ and $\lieg_{1}$ are the subalgebras with nonzero entries only in the lower-left and upper-right blocks, respectively.
The decomposition 
satisfies $[\lieg_i,\lieg_j]\subset\lieg_{i+j}$, where we set $\lieg_k=\{0\}$ for
$|k|>1$.
Note that $\lieg_{-1} \cong \Hom(\BR^m, \BR^n)$, and the restriction of
the adjoint action of elements of $G_0$ to $\lieg_{-1}$ is the representation $\rho$ from Section \ref{AGr}.  


Next, let $P < G$ comprise the
block-upper-triangular matrices,
with Lie algebra $\liep:=\lieg_0\ltimes\lieg_1\subset\lieg$. It is the
stabilizer of $\BR^m \subset \BR^{m+n}$, 
so $G/P$ can be identified with $\Gr(m,n)$.
As is well known, the tangent bundle
$T(G/P)$ is the associated bundle $G \times_P (\lieg / \liep)$, where $P$ acts via the adjoint representation, which factors on the quotient $\lieg / \liep$ through projection to $G_0$.
The vector space $\lieg/\liep$ is moreover $G_0$-equivariantly isomorphic to $\lieg_{-1}$. 
Consequently, $\Gr(m,n)$ carries an 
almost
Grassmannian structure.
 Note that the
auxiliarly bundles $\Cal E \cong G \times_P \BR^m$ and $\Cal F \cong G \times_P (\BR^{m+n}/\BR^m)$ for this structure are exactly the tautological and
the anti-tautological bundles. 

\subsection{Automorphisms and flatness}\label{sec.auts}

Let $M$ be a $C^{k+1}$ manifold of dimension $mn$ with a $C^k$ almost-Grassmannian structure comprising
\begin{itemize}
\item $C^{k+1}$ vector bundles $\mathcal{E}$ and $\mathcal{F}$ of ranks $m$ and $n$, respectively
\item $\Theta : TM \stackrel{\sim}{\rightarrow} \mathcal{E}^* \otimes \mathcal{F}$ of regularity $C^k$
\item $\nu : \La^m \mathcal{E} \stackrel{\sim}{\rightarrow} \La^n \mathcal{F}$ of regularity $C^k$
\end{itemize}

\begin{definition} \label{def-aut}
An \emph{automorphism} of the $C^k$ almost-Grassmannian structure above is a $C^{k+1}$ diffeomorphism $h$ of $M$ together with
\begin{itemize}
\item lifts $h_{\mathcal{E}}^*$ and $h_{\mathcal{F}}$ of $h$ to automorphisms of $\mathcal{E}^*$ and $\mathcal{F}$, respectively, 
\item such that $h^* \Theta = (h_{\mathcal{E}}^* \otimes h_{\mathcal{F}}) \circ \Theta$, and
\item such that $\nu \circ \La^m h_{\mathcal{E}} = \La^n h_{\mathcal{F}} \circ \nu$ 
\end{itemize}
\end{definition}

In the $G$-structure framework, $h \in \mbox{Diff}^{k+1} M$ is an automorphism if it lifts to a principal bundle automorphism of $\Cal G_0$ which is semi-conjugate via the $\rho$-equivariant bundle morphism $\Cal G_0 \to \mathcal{P}M$ to the natural lift of $h$ to $\mathcal{P}M$.  

Isomorphisms of almost-Grassmannian structures are defined by the obvious extension of Definition \ref{def-aut}.  Local isomorphisms are isomorphisms between connected open subsets, with their restricted structures.

Consider $g \in G$ as a diffeomorphism of $\Gr(m,n) \cong G/P$.  It naturally acts by automorphisms $g_{\mathcal{E}}^*$ and $g_{\mathcal{F}}$ of the vector bundles $\mathcal{E}^* \cong G \times_P \BR^{m*}$ and ${\mathcal{F}} \cong G \times_P (\BR^{m+n}/\BR^m)$, respectively.  The $P$-equivariant isomorphism of $\Hom(\BR^m,\BR^{m+n}/\BR^m)$ with $\lieg / \liep$ gives a $G$-equivariant isomorphism of $T \Gr(m,n) \cong G \times_P \lieg/\liep$ with  $\mathcal{E}^* \otimes \mathcal{F}$, on which $g$ acts by $g_{\mathcal{E}}^* \otimes g_{\mathcal{F}}$.  Any $g \in G$ is thus an automorphism of the almost-Grassmannian structure on $\Gr(m,n)$.

Let $P_+$ be the connected, unipotent, normal subgroup of $P$ with Lie algebra $\lieg_{1}$.  For $g \in P_+$, the linear isomorphisms $(g_{\mathcal{E}}^*)_{[\Id_G]}$ and $(g_{\mathcal{F}})_{[\Id_G]}$ are trivial, so the derivative of $g$ on $T_{[\Id_G]} (G/P)$ is trivial.  The $G$-conjugates of $P_+$ furnish nontrivial strongly essential automorphisms at every point of $\Gr(m,n)$.

Each $g \in G \backslash \{ \pm \Id_{\BR^{m+n}} \}$ is a nontrivial transformation of $G/P$, because there is no larger $G$-normal subgroup in $P$.  Thanks to the canonical Cartan connection associated to an almost-Grassmannian
structure, presented in Section \ref{harm-curv} below, the automorphism group of $\Gr(m,n)$ is precisely $G/\{ \pm \Id \}$.  On any almost-Grassmannian manifold (of sufficient regularity),
the Cartan connection underlies the fact that the automorphism group is a Lie group of dimension at most $(m+n)^2-1 = \dim G$, with equality only if the
structure is locally isomorphic to $\Gr(m,n)$, in which case it is said to be \emph{flat}. Deciding whether an almost
Grassmannian structure is flat thus is a fundamental question
in the theory.

\subsection{The harmonic torsion}\label{tors}
The description as a $G_0$-structure $(p_0:\Cal G_0\to M,\theta)$ directly leads to the first fundamental
invariants of almost Grassmannian structures. 
We first choose a $C^{k}$ principal connection
$\gamma\in\Omega^1(\Cal G_0,\lieg_0)$. If $\theta$ is at least $C^1$, we can define the \textit{torsion} of $\gamma$ as the covariant exterior
derivative $d^\gamma\theta\in\Omega^2(\Cal G_0,\BR^{mn})$; explicitly, for $\xi,\eta\in\frak X(\Cal G_0)$,
\begin{equation}
\label{eqn.dtheta} 
d^\gamma\theta(\xi,\eta)=d\theta(\xi,\eta)+\gamma(\xi)(\theta(\eta))-\gamma(\eta)(\theta(\xi)).
\end{equation}
If $\theta$, $\xi$ and $\eta$ are $C^k$, then the above is a $C^{k-1}$ function.  
From the properties of $\theta$, it follows readily that $d^\gamma\theta$ is horizontal and
$G_0$-equivariant and thus descends to a form $T^\gamma\in\Omega^2(M,TM)$, which is
the usual interpretation of the torsion. 

We next compute the dependence of $T^\gamma$ on $\gamma$.  First note that, at a point $u \in \Cal G_0$, (\ref{eqn.dtheta}) depends only on $\gamma_u:T_u\Cal G_0\to\lieg_0$.   As discussed in \ref{model}, we can view $\theta$ as having values in $\lieg_{-1}$.  For any other principal connection $\hat{\gamma}$, the 
difference $\hat\gamma_u - \gamma_u$ is given by $f_u \circ\theta_u$ for some linear map
$f_u :\lieg_{-1}\to\lieg_0$.  The first differential in the cochain complex of $\lieg_{-1}$ with coefficients in $\lieg$ restricts on $\lieg_{-1}^* \otimes \lieg_0$ to the following $G_0$-equivariant linear map:
\begin{equation}
  \label{partial1}
\partial_1:\lieg_{-1}^*\otimes\lieg_0\to\La^2\lieg_{-1}^*\otimes\lieg_{-1}   \qquad (\partial_1f)(w,v) = f(w) v - f(v) w
\end{equation}
For all $u$,
\begin{equation}
\label{eqn.change.torsion}
T^{\hat\gamma}_u - T^\gamma_u =  (\partial_1
f_u) \circ \theta_u.
\end{equation}
  The image of $\partial_1$ 
determines a
smooth subbundle $\Cal A\subset\La^2T^*M\otimes TM$. The projection of $T^\gamma$ to
$(\La^2T^*M\otimes TM)/\Cal A$ is thus independent of the
choice of connection (see \cite[Secs 3.1.10--3.1.13]{cap.slovak.book.vol1}). This invariant of the almost Grassmannian structure is called the \textit{intrinsic torsion} or the \textit{harmonic torsion}.

For $\Gr(m,n)$ it is easy to see that locally
there always are torsion-free connections preserving the structure, so the intrinsic
torsion of the homogeneous model vanishes identically. Thus nonzero intrinsic torsion is an obstruction to local isomorphism of a given
almost Grassmannian structure to $\Gr(m,n)$. 
For a $C^1$-structure, it is an
obstruction to local $C^1$-isomorphism to $\Gr(m,n)$ (for which the
corresponding map between the underlying manifolds would be a local $C^2$-diffeomorphism).

We now explicitly describe the subbundle $\Cal A\subset\La^2T^*M\otimes TM$ as a $\lieg_0$ representation when $m=2$.
Recall that
the representation corresponding to the tangent bundle
$TM$ is $\lieg_{-1} \cong \BR^{2*}\boxtimes\BR^n$, where the exterior tensor product corresponds to the direct sum decomposition of $\lieg_0$. 
Next we have the decomposition into irreducible
components
\begin{equation}\label{La2-decomp}
  \La^2(\lieg_{-1}^*)\cong (\La^2\BR^2\boxtimes S^2\BR^{n*})\oplus
  (S^2\BR^2\boxtimes\La^2\BR^{n*}).
\end{equation}
We tensor these with $\lieg_{-1}$ and decompose into irreducibles. 
For $k\geq 2$,
the representation $S^2\BR^k\otimes\BR^{k*}$ splits into a trace-free component, denoted
$(S^2\BR^k\otimes\BR^{k*})_0$, and a trace component, isomorphic to $\BR^k$, and similarly
for the dual. There is an analogous decomposition of $\La^2\BR^k\otimes\BR^{k^*}$,
but here the trace-free part is trivial when $k=2$.

For $m=n=2$, the map $\partial_1$ from \eqref{partial1} is
surjective, so 
no intrinsic torsion is
available. In our case when $m=2$, $n>2$,
$$(\La^2\BR^2\otimes\BR^{2*})\boxtimes (S^2\BR^{n*}\otimes\BR^n)\subset\Im(\partial_1).$$
The intersection of $\Im(\partial_1)$ with the other irreducible components of $\La^2(\lieg_{-1}^*) \otimes \lieg_{-1}$ is the trace component, which can be written
\begin{equation}
  \label{tracep}
\BR^2\boxtimes (\La^2\BR^{n*}\otimes\BR^n)+
(S^2\BR^2\otimes\BR^{2*})\boxtimes\BR^{n*},  
\end{equation}
where the factors $\BR^2$ and $\BR^{n*}$ are embedded via a tensor product with $\Id$ followed by a symmetrization and an alternation, respectively. Hence 
the harmonic torsion 
corresponds to a section of the bundle associated to the remaining irreducible component
$$
\mathbb{T} = (S^2\BR^2\otimes\BR^{2*})_0\boxtimes (\La^2\BR^{n^*}\otimes\BR^n)_0.
$$
To verify non-vanishing harmonic torsion in the example we are going to construct, we
need the following result.
\begin{lemma}\label{harm-tors}
  Let $\xi,\eta\in \Hom(\BR^2,\BR^n)$ both have 
  kernel spanned by $0 \neq v\in\BR^2$, and let $T \in \Im(\partial_1)$. 
 Then
  $T(\xi,\eta):\BR^2\to\BR^n$ maps $v$ into the span of the images of
  $\xi$ and $\eta$.
\end{lemma}
\begin{Pf}
  Take $0 \neq \al\in\BR^{2*}$ with $\al(v)=0$, so we
  can write $\xi=\al\otimes w_1$ and $\eta=\al\otimes w_2$ for elements
  $w_1,w_2\in\BR^n$. 
  Let $T=T_1+T_2$ be the decomposition of $T$ corresponding to \eqref{La2-decomp}. 
  Clearly $T_1(\xi,\eta)=0$. 
 Now decompose $T_2$ according to 
  \eqref{tracep}
  as $T_{21}+T_{22}$ (in a non-unique way).  Given $\tilde v\in\BR^2$, embedded in the trace component of $S^2\BR^2\otimes\BR^{2*} \cong \Hom(\BR^2,S^2\BR^2)$, it sends $v$ to a multiple
  of the symmetric product $v\odot\tilde v$.  Since $\xi(v)=\eta(v)=0$, we conclude that
  $T_{1}(\xi,\eta)(v)=0$. On the other hand, $\tau\in\BR^{n^*}$, embedded in the trace component of $\La^2\BR^{n*}\otimes\BR^n \cong \Hom(\La^2 \BR^n, \BR^n)$,
  sends $(w_1,w_2)$ to a
  multiple of $\tau(w_1)w_2-\tau(w_2)w_1$, so all values of $T_{2}(\xi,\eta)$ belong to the span of the images of $\xi$ and $\eta$.  The desired conclusion follows.
\end{Pf}

\subsection{Deformations of almost-Grassmannian structures}\label{sec.inf.def}

Given an al\-most-Grassmannian structure with $\theta : TM \stackrel{\sim}{\to} \mathcal{E}^* \otimes \mathcal{F}$, we will construct deformations by post-composing with a continuous family $\{ \Phi_t \}$ of linear automorphisms of $\mathcal{E^*} \otimes \mathcal{F}$.  To construct this family, we will first construct endomorphisms, that is, a section $\Phi$ of $\End \mathcal{E}^* \otimes \End \mathcal{F}$ and then show that this exponentiates to a one-parameter family of automorphisms.

At a given point $x \in M$, write $\mathcal{E}_x \cong E$ and $\mathcal{F}_x \cong F$.  The vector space automorphisms of $E^* \otimes F$ respecting the tensor product are those of the form $\Psi_{E^*} \otimes \Psi_F$, for $\Psi_{E^*} \in \Aut E^*$ and $\Psi_F \in \Aut F$.  Given a one-parameter group of such automorphisms $\Psi^t_{E^*} \otimes \Psi^t_F$, the generating endomorphism has the form $\psi_{E^*} \otimes \Id_F + \Id_{E*} \otimes \psi_F$ for $\psi_{E*} \in \End E^*$ and $\psi_F \in \End F$.  The condition $\det \Psi^t_{E^*} \cdot \det \Psi^t_F \equiv 1$ is equivalent to $\tr\  \psi_{E*} + \tr\  \psi_F =0$.

The sections of $\mbox{Aut}(\mathcal{E}^* \otimes \mathcal{F})$ arising from automorphisms of the almost-Grassmannian structure are those of the form $\Psi_{\Cal E^*} \otimes \Psi_{\Cal F}$, for $\Psi_{\Cal E^*} \in \Aut \Cal E^*$ and $\Psi_{\Cal F} \in \Aut {\Cal F}$, with $\nu \circ \La^m \Psi_{\Cal E^*} = \La^n \Psi_{\Cal F} \circ \nu$.  The generator of a nontrivial deformation is thus nontrivial modulo $\End \Cal E^* \otimes \Id_{\Cal F} + \Id_{\Cal E^*} \otimes \End \Cal F$.  A pointwise complementary subbundle is given by the tensor product of trace-free endomorphisms $\mbox{End}_0 \  \Cal E^* \otimes \mbox{End}_0 \  \Cal F$.  We will construct a section of this bundle in Sections \ref{sec.eigen.sections} and \ref{sec.endo.section} below.

The results of \cite{cap.deformations} apply to almost-Grassmannian structures of sufficient regularity to define a Cartan connection (see Section \ref{harm-curv}).  Here infinitesimal automorphisms and deformations are described as the kernel and cokernel, respectively, of BGG operators acting on sections of the adjoint tractor bundle, with a certain ``twisted'' linear connection.  The infinitesimal change of harmonic torsion and harmonic curvature (for the latter, see Section \ref{harm-curv}) resulting from a given infinitesimal deformation can also be described in general from these operators and this connection.  This point of view was the inspiration for the concrete deformations we construct below.

\subsection{Prolongation and the canonical Cartan connection}\label{harm-curv}

We will verify in Section \ref{sec.C2.curvature} that the results of \cite{cap.me.parabolictrans} apply to $C^k$ almost-Grassmannian structures with $k \geq 2$, 
so any example of this regularity admitting a strongly essential flow by automorphisms
has vanishing harmonic curvature on an open set containing the higher-order fixed point in its closure.  We explain in this section that $(2,n)$-almost Grassmanian structures of regularity $C^k$ with $k\geq 2$ determine a canonical $C^0$ Cartan connection as well as $C^0$ harmonic curvature.  In low regularity, general existence results do not apply, so we briefly sketch the explicit constructions, following \cite{css.ahs2}.

\subsubsection{Construction of the prolongation}

Given a $C^k$ almost Grassmannian structure $(p_0:\Cal G_0\to M,\theta)$ of type $(2,n)$, we will prolong $\Cal G_0$ to a $C^{k-1}$ principal $P$-bundle $\mathcal{G} \to M$.  
To this end, we view $\lieg_0$ as a subalgebra of $\End \lieg_{-1}$.
The kernel $\ker(\partial_1)$ of the differential from (\ref{partial1}) is the subspace of $\lieg_{-1}^*\otimes\lieg_0 \subset \lieg_{-1}^* \otimes \lieg_{-1}^* \otimes \lieg_{-1}$ of elements symmetric in $\lieg_{-1}^*$, which is precisely the \emph{first prolongation} of $\lieg_0$  (see \cite[I.1]{kobayashi.transf}). It is isomorphic to
$\lieg_1$, embedded into $\lieg_{-1}^*\otimes\lieg_0$ via the adjoint representation. 

The bundle $\Cal G$ is constructed as a $\lieg_1$-bundle over $\Cal G_0$. Note that $P \cong G_0 \ltimes \lieg_1$.
Given $u_0\in\Cal G_0$,
denote by $V_{u_0} \Cal G_0\subset T_{u_0}\Cal G_0$
the vertical subspace, so $\theta_{u_0}$ defines a linear isomorphism
$T_{u_0}\Cal G_0/V_{u_0}\Cal G_0\to\BR^{2n}$.  Recall from Section \ref{tors} that (\ref{eqn.dtheta}) depends at $u_0$ only on the value $\gamma_{u_0}$ of a chosen principal connection. Let $u : T_{u_0}\Cal G_0 \to\lieg_0$ be any linear
map recognizing fundamental vector fields---that is $u(\zeta_A(u_0)) = A\in\lieg_0$ for $\zeta_A(u_0) = \left. \frac{\mathrm{d}}{\mathrm{d}t}\right|_0 u_0. e^{tA}$. Now $d \theta_{u_0} + [u,\theta_{u_0}]$ vanishes when either input is in $V_{u_0}\Cal G_0$, 
so it equals $\theta_{u_0}^* T_u$ for
a unique map
$T_u\in\La^2\lieg_{-1}^*\otimes\lieg_{-1}$. 

Varying the choice of $u$ yields, as in Section \ref{tors}, 
the affine subspace $T_u+\Im(\partial_1)$;
moreover, the trace-free subspace
$\mathbb{T} = (S^2\BR^{2*}\otimes\BR^2)_0\boxtimes (\La^2\BR^n\otimes\BR^{n*})_0$ is a $G_0$-invariant
complement to $\Im(\partial_1)$.  Thus for each point $u_0\in\Cal G_0$, the linear map $u$ can be chosen
such that $T_u$ is totally trace-free.  For a fixed $u_0$ and $T \in \mathbb{T}$, the space of linear maps $u$ giving rise via $\theta_{u_0}$ to $T$ is, according to (\ref{eqn.change.torsion}), an affine space modeled on $\ker \partial_1 \cong \lieg_{1}$.
Explicitly,
any two such maps differ according to $\hat u - u = \mbox{ad}_Z \circ \theta_{u_0} \in \lieg_0$, for a unique
element $Z\in\lieg_1$.

Now $\Cal G$ is the family of maps $u$ as above for which $T_u \in \mathbb{T}$, as $u_0 $ varies over $\Cal G_0$.  Denote $q$ the projection $\Cal G \to \Cal G_0$.
We can realize $\Cal G$ as a subspace of the vector bundle $T^*\Cal G_0\otimes\lieg_0\to\Cal G_0$ 
defined as above in terms of $\theta$ and $d\theta$, which is $C^{k-1}$;
 it follows that $\Cal G$ is a $C^{k-1}$ submanifold here.
Elementary representation theory gives a $G_0$-equivariant linear map
$S:\La^2\lieg_{-1}^*\otimes\lieg_{-1}\to\lieg_{-1}^*\otimes\lieg_0$ which vanishes on
$\mathbb{T}$ and such that $\partial_1\o S$ is the projection to
$\Im(\partial_1)$.  From a $C^{k-1}$ principal connection on $\Cal G_0$, one can use $S$ to modify it to a $C^{k-1}$ principal connection
$\gamma$ with $\gamma_{u_0} \in \Cal G$ for all $u_0 \in \Cal G_0$; connections of this type correspond to local $C^{k-1}$ sections of $q$. 
Now $q:\Cal G\to\Cal G_0$ is a $C^{k-1}$ principal $P_+$-bundle, and $p:=p_0\o q:\Cal G\to M$ is  a $C^{k-1}$
principal $P$-bundle
(see \cite{css.ahs2} for more details).

\subsubsection{Harmonic curvature and Cartan connection}
\label{sec.curv.and.connxn}

Construction of the Cartan connection on $\Cal G$ entails, by analogy with the prolongation process of the previous section, finding canonical $\lieg_1$-valued one-forms on $\Cal G$, which turn out to be unique.  There are tautological forms $\theta_{-1} + \theta_0$, where $\theta_{-1}:=q^*\theta$ and $(\theta_0)_u:= q^* u$, viewing $u \in \Hom(T_{q(u)} \Cal G_0,\lieg_0)$.
It is easy to see that
$(\theta_0)_u(\zeta_A)=A$ for all $A\in\lieg_0$,
and that $\theta_{-1} + \theta_0$ is $P$-equivariant once
 $\lieg_{-1}\oplus\lieg_0$ is identified with $\lieg/\lieg_1$.

Assuming that $\theta$ is at least $C^2$, the form $\theta_0$ is at least
$C^1$, so we can form its exterior derivative $d\theta_0 \in \Omega^2(\Cal G, \lieg_0)$.
Let 
$\phi:T_u\Cal G\to\lieg_1$ be a linear map satisfying
$\phi(\zeta_{A+Z}(u))=Z$ on the fundamental vector fields for $A+Z\in  \lieg_0\ltimes\lieg_1 \cong \liep $.  
As before, 
\begin{equation}\label{K_u}
(d \theta_0)_u + \frac{1}{2} [ (\theta_0)_u, (\theta_0)_u ] + [\phi,(\theta_{-1})_u]
\end{equation}
vanishes if either input is in $V_u \Cal G$ (the vertical bundle for $p : \Cal G \to M$), so 
it equals $(\theta_{-1})^* K_\phi$ for a linear map
$K_\phi:\Lambda^2\lieg_{-1}\to\lieg_0$.  

For another choice $\hat{\phi}$, the difference is
$\hat\phi - \phi = f \circ (\theta_{-1})_u$ for some linear map
$f\in\lieg_{-1}^*\otimes\lieg_1$.  Another grading component of our Lie algebra differential,
$\partial_2:\lieg_{-1}^*\otimes\lieg_1\to\Lambda^2\lieg_{-1}^*\otimes\lieg_0$,  is a $G_0$--equivariant linear map for which $K_{\hat\phi}-K_\phi=\partial_2(f)$. 
Projection modulo the subbundle $\mathcal{B} \subset \Lambda^2T^*\Cal G\otimes\lieg_0$ determined by $\Im(\partial_2)$ yields another invariant: 
given a local principal connection $\gamma$ on $\Cal G_0$ for which $\gamma_{u_0} \in \Cal G$ for all $u_0$ in the domain, and a locally defined $\phi \in \Omega^1(\Cal G, \lieg_1)$,  smooth of class $C^{k-1}$,
one obtains a $C^{k-1}$ section of $(\Lambda^2T^*\Cal G\otimes\lieg_0) / \mathcal{B}$ called the \textit{harmonic curvature} of the geometry. This can be
recovered as a component of the curvature of any adapted connection.

Now, it can be shown that $\partial_2$ is injective, and that there is a
natural $G_0$--invariant complement $\mathbb{K}$ to $\Im(\partial_2)$. Hence for each $u\in\Cal G$, there
is a unique $\phi$ such that $K_\phi \in \mathbb{K}$. We obtain $\theta_1\in\Omega^1(\Cal G,\lieg_1)$ of class $C^{k-2}$, and a $C^{k-2}$ Cartan connection 
$\omega=\theta_{-1}\oplus\theta_0\oplus\theta_1 \in \Omega^1(\Cal G, \lieg)$.
Note that when $k \geq 3$, the Cartan curvature can be defined by
$K=d\omega+ \frac{1}{2} [\omega,\omega] \in \Omega^2(\Cal G,\lieg)$,
and the harmonic torsion and the harmonic curvature
are components of $K$. For the homogeneous model, $\omega$ is the Maurer-Cartan form, so $K$ vanishes identically, as does the harmonic
curvature (see \cite{css.ahs2} for more details).

\begin{proposition}
\label{prop.cg.kgeq2}
For $k \geq 2$, a $C^k$ almost-Grassmannian structure of type $(2,n)$, $n \geq 3$, on $M$ determines a $C^{k-1}$ principal $P$-bundle $\Cal G \to M$ equipped with a $C^{k-2}$ Cartan connection $\omega \in \Omega^1(\Cal G, \lieg)$.  A $C^k$ morphism between two such almost-Grassmannian structures canonically lifts to a morphism of the associated Cartan geometries.
\end{proposition}


\begin{Pf}
It remains only to verify the last statement.
Let $h$ be a local $C^{k+1}$ diffeomorphism between open subsets of $M$ and $\widetilde{M}$, lifting to 
 a $G_0$-equivariant local $C^k$-diffeomorphism 
 $\Phi_0:\Cal G_0\to\widetilde{\Cal G}_0$ with
$\Phi^*_0 \tilde\theta=\theta$.
Given $u\in q^{-1}(u_0) \subset \Cal G$, 
let
$\Phi(u):= (\Phi_0^{-1})^* u$.  It is easy to check that $T_{\Phi(u)}=T_u$,
so
$\Phi(u)\in (\tilde{q})^{-1}(\Phi_0(u_0)) \subset  \tilde{\Cal G}$. 
This
evidently defines a $P$-equivariant $C^{k-1}$ lift $\Phi: \Cal G\to\tilde{\Cal G}$ of $\Phi_0$. 
The construction implies that $\Phi^*\tilde\theta_i=\theta_i$
for $i=-1,0$. 
Then repeating this argument with the harmonic curvature allows us
to conclude that $\Phi^*\tilde\omega=\omega$, so
$\Phi$ is a morphism of Cartan geometries. 
\end{Pf}

As a corollary, we note that for a
structure of class at least $C^2$, nonvanishing harmonic curvature is an obstruction
to local isomorphism to $\Gr(2,n)$.

\section{Description of the strongly essential flow in coordinates}
\label{sec.descr.flow}

By homogeneity of $\Gr(2,n)$, we may assume the point $x_0$ in Theorem \ref{main.thm} is the standard 2-plane spanned by the first two coordinate vectors in $\BR^{2+n}$.
The deformation will be constructed on the open subset $U$ comprising the orbit of $x_0$ under all transformations $\Id + X$, with $X \in \mbox{Hom}(\BR^2, \BR^n)$.  This set is the domain of an affine chart in the Pl\"ucker coordinates on $\Gr(2,n)$, and equals the top cell in the standard Schubert decomposition.  Identify $U$ with $\Hom(\BR^2, \BR^n)$, and represent an element of $U$ in coordinates 
$$X = ( x_{ij})_{i=1,j=1}^{i=n,j=2}$$

The principal $P$-bundle $G \rightarrow \Gr(2,n) \cong G/P$ restricted to $U$ is smoothly equivalent to the trivial bundle $U \times P$.  The quotient bundle $G/G_1$ restricted to $U$ is $\left. \Cal G_0 \right|_U$, which is smoothly equivalent to $U \times G_0$.

 We will use the following explicit trivializations over $U$ of the tautological and anti-tautological bundles, together with their duals.
 Denote by  $\{e_1, \ldots, e_{n+2} \}$ the standard basis of $\BR^{2+n}$ with dual basis $\{e^1, \ldots, e^{n+2} \}$.  
 For $j'=1',2'$, define a section of  $\left. \mathcal{E}\right|_U$ by $E_{j'}(X) = (\Id + X)e_{j'}$, and let $E^{j'}$ be the sections $e^{j'}$ of $\mathcal{E}^*$.  Next let $E_i(X)$ equal the image of $e_{i+2}$ in $\BR^{2+n}/X$ for $i=1,\ldots, n$, and $E^i(X) = e^{i+2}  - x_{i1} e^1 - x_{i2} e^2$, which are well-defined on $\BR^{2+n} / X$.
 We will henceforth denote the restrictions of these bundles to $U$ simply by $\mathcal{E}$, $\mathcal{E}^*$, $\mathcal{F}$, and $\mathcal{F}^*$.


The restriction of the standard flat Grassmannian structure on $\Gr(2,n)$ to $U$ is given by the obvious isomorphism from $TU\cong U\x\Hom(\BR^2, \BR^n)$ with $\Cal E^* \otimes \Cal F$. It sends the coordinate vector fields $\partial^{j'}_{i}:= {\partial}/{\partial x_{ij'}}$ on $U$ to the sections $E^{j'}\otimes E_i$ of $\Cal E^* \otimes \Cal F$. 

\subsection{The strongly essential flow}

Let $Z$ be a rank-one element of $\lieg_1 \cong \Hom(\BR^n, \BR^2)$.  Let $\{ z^t \}$ be the one-parameter subgroup of $P <  G$ generated by $Z$; it is just the group of matrices $\{ \Id + tZ \} < \SL(2+n,\BR)$.  Theorem \ref{main.thm} applies to any strongly essential flow generated by a rank-one element of $\lieg_{1}$.  After conjugation in $G$, we may assume $Z = e^1 \otimes e_{1'}$, where now $\{ e^1, \ldots, e^n \}$ is the standard basis of $\BR^{n*}$ and $\{ e_{1'}, e_{2'} \}$ the standard basis of $\BR^2$.  Then $\im Z = \BR e_{1'}$, and $\ker Z = \mbox{span}\{ e_2, \ldots, e_n \}$.

For $X \in U \cong \lieg_{-1}$, denote $e^X$ the corresponding lower-triangular unipotent matrix in $G$.  Then compute the image in $U \times P$
\begin{eqnarray*}
z^t.e^X & = & 
\left( 
\begin{array}{cc} 
\Id_2 + tZX & tZ \\
X & \Id_n
\end{array} 
\right)  \\
  & = &  
\left(
\begin{array}{cc}
\Id_2 & 0 \\
X(\Id_2 + tZX)^{-1} & \Id_n 
\end{array}
\right) 
\left(
\begin{array}{cc}
\Id_2 + tZX & tZ \\
0 & \Id_n - X (\Id_2 + tZX)^{-1} tZ
\end{array}
\right),
\end{eqnarray*}

assuming that $\Id_2 + tZX$ is invertible, which for fixed $t$ holds on an open neighborhood of $0$.

The following two subspaces are fixed by $\{z^t\}$: 
$$ F_1 = \{ X \ : \ XZ = 0 \}$$
and 
$$ F_2 = \{ X \ : \ ZX = 0 \}$$
The intersection $F_1 \cap F_2$ will be called the \emph{strongly fixed set} and denoted $SF$.  Note that $X \in SF$ if and only if $[X,Z]=0$ and, in coordinates,
$$SF = \{ X \ : \ x_{12} = 0 = x_{i1} \ \mbox{for all} \  i = 1, \ldots, n \}$$


Let $H_0 = \{ X \ : \ x_{11} = 0 \}$.  For $X \notin H_0$, decompose $X$ as
$$X_f + X_d = 
\left(
\begin{array}{cc}
0 & 0 \\
    & x_{22} - \frac{x_{12}}{x_{11}} \cdot x_{21} \\
 \vdots & \vdots \\
 0 & x_{n2} - \frac{x_{12}}{x_{11}} \cdot x_{n1}
 \end{array}
\right) 
+
\left(
\begin{array}{cc}
x_{11} & \frac{x_{12}}{x_{11}} \cdot x_{11}  \\
  &  \\
 \vdots & \vdots \\
 x_{n1} & \frac{x_{12}}{x_{11}} \cdot x_{n1}
 \end{array}
\right) 
$$

with $X_f \in SF$ and $X_d$ rank 1.
If $X \notin H_0$, then $z^t.X = X_f + z^t.X_d$, which equals
\begin{eqnarray}
\label{eqn.ztaction}
X_f + 
\left(
\begin{array}{cc}
\frac{x_{11}}{1+tx_{11}} & \frac{\frac{x_{12}}{x_{11}} \cdot x_{11}}{1+tx_{11}} \\
\vdots & \vdots \\
\frac{x_{n1}}{1+tx_{11}} & \frac{\frac{x_{12}}{x_{11}} \cdot x_{n1}}{1+tx_{11}}
\end{array}
\right)
\end{eqnarray}

Let $H_+ = \{ X \in U \ : \ x_{11} > 0 \}$ and $H_- = \{ X \in U \ : \ x_{11} < 0 \}$.  The formula (\ref{eqn.ztaction}) above yields $z^t.X \rightarrow X_f$ as $t \rightarrow \pm \infty$ for $X \in H_{\pm}$, respectively. 

Note that if $X \in H_0$, then $ZXZ = 0$, and $(\Id + tZX)^{-1} = \Id - tZX$.  Then $z^t X = X (\Id - t ZX) $ which equals $X$ as $t$ varies if and only if $XZX=0$; the latter holds only when $ZX$ or $XZ=0$.  We conclude that $F_1 \cup F_2$ equals the fixed set of ${z^t}$ in $U$.  

\subsection{Action on associated vector bundles}
\label{sec.assoc.vbs}

The matrix in $P$
$$ p_t(X) = 
\left(
\begin{array}{cc}
\Id_2 + tZX & tZ \\
0 & \Id_n - X (\Id_2 + tZX)^{-1} tZ
\end{array}
\right)
$$ 
from above encodes the action of $z^t$ on $\mathcal{E}$ and $\mathcal{F}$, and, in turn, on $TU$.  For $X \notin H_0$, 
 \begin{eqnarray}
\label{eqn.holonomy}
p_t(X) = p_t(X_d) =
\left(
\begin{array}{cccccc}
1+tx_{11} & tx_{12} & t                                        &  &  \cdots & 0  \\
0               & 1           & 0                                        &   &  \cdots & 0 \\
                 &              & \frac{1}{1+tx_{11}}          &   &  &   \\
                 &              & \frac{-tx_{21}}{1+tx_{11}} & 1 & &   \\
                 &              &      \vdots                            &     & \ddots &     \\
                 &               & \frac{-tx_{n1}}{1+tx_{11}}&  &               & 1
\end{array}
\right)                
\end{eqnarray}

For $X \in H_0$, straightforward calculation gives the formula (\ref{eqn.holonomy}), with $x_{11} = 0$.

On $\mathcal{E} = \left. (G \times_P \BR^2) \right|_U \cong U \times \BR^2$ the action of $\{ z^t \}$ is
\begin{eqnarray}
\label{eqn.E.action}
(z^t_{\mathcal{E}})_X = 
\left(
\begin{array}{cc}
1+tx_{11} & tx_{12} \\
0 & 1 
\end{array}
\right) \ \mbox{with respect to} \ \{ E_{1'}, E_{2'} \}
\end{eqnarray} 

and on $\mathcal{E}^*$,
\begin{eqnarray}
\label{eqn.Estar.action}
(z^{-t}_{\mathcal{E}})_{z^t.X}^* = 
\left(
\begin{array}{cc}
\frac{1}{1+tx_{11}} & \frac{-tx_{12}}{1+tx_{11}} \\
0 & 1 
\end{array}
\right) \ \mbox{with respect to} \ \{ E^{1'}, E^{2'} \}
\end{eqnarray}

On $\mathcal{F}$, the action is
\begin{eqnarray}
\label{eqn.F.action}
(z^t_{\mathcal{F}})_{X} = 
\left(
\begin{array}{cccc}
\frac{1}{1+tx_{11}}          &   &  &   \\
\frac{-tx_{21}}{1+tx_{11}} & 1 & &   \\
   \vdots                            &     & \ddots &     \\
 \frac{-tx_{n1}}{1+tx_{11}}&  &               & 1
\end{array}
\right) \ \mbox{w.r.t.} \ \{ E_1, \ldots, E_n\}
\end{eqnarray}

and, finally, on $\mathcal{F}^*$, 
\begin{eqnarray}
\label{eqn.Fstar.action}
(z^{-t}_{\mathcal{F}})^*_{z^t.X} = 
\left(
\begin{array}{cccc}
1+tx_{11} &  &  &  \\
tx_{21} & 1 &   &   \\
\vdots &    & \ddots & \\
tx_{n1} &   &      & 1
\end{array}
\right) \ \mbox{w.r.t.} \ \{ E^1 , \ldots, E^n \}
\end{eqnarray}

\section{The invariant deformation and non-flatness}

Our deformation is constructed in Sections \ref{sec.eigen.sections} and \ref{sec.endo.section} below.  Then in Section \ref{sec.proof.calc} we compute terms of the harmonic torsion which are nonzero on an open, dense subset of any neighborhood of $x_0$.  It follows that our $C^1$, deformed structures are not flat on any open set $V$ with $x_0 \in \overline{V}$, and thus that the $C^1$ version of \cite[Thm 1.3]{mn.1graded} does not hold.  The proof of vanishing harmonic torsion in \cite{mn.1graded} requires several degrees of differentiability of the structure, so it remains open whether $C^1$ is the maximal regularity of such a counterexample.  
Our result \cite[Prop 3.5]{cap.me.parabolictrans}, on the other hand, says that the harmonic curvature must vanish on an open set $V$ with $x_0 \in \overline{V}$, in the presence of a flow by strongly essential autmorphisms, and we explain in Section \ref{sec.C2.curvature} below that it applies to $C^2$ structures.  In Section \ref{sec.nonzero.curv}, we verify that the harmonic curvature of our deformations restricted to their common smooth locus is nonzero.  Our deformations are thus in some sense structures of maximal regularity for which the conclusion of \cite[Prop 3.5]{cap.me.parabolictrans} does not hold.

\subsection{Eigen-sections of associated bundles}
\label{sec.eigen.sections}

Here we define sections of $\mathcal{E}$ and $\mathcal{F}$ and of the dual bundles, which are invariant by $z^t$ up to multiplication by a function on $U$.  For any $X$ where the decompositions
\begin{equation}
\label{eqn.decomps.EF}
\mathcal{E}_X = \ker X_d \oplus \im Z  \qquad \mbox{and} \qquad \mathcal{F}_X = \ker Z \oplus \im X_d 
\end{equation}
are valid, each section will have values in one factor or its dual.  The sections together will span the fibers over $X$ in $\mathcal{E}, \mathcal{F}$, or their duals.  

Define
\begin{eqnarray*}
v(X) = - x_{12} E_{1'} + x_{11} E_{2'}  &  & \iota(X) = E_{1'} \\
\tilde{v}(X) = E^{2'} &  & \tilde{\iota}(X) = x_{11} E^{1'} + x_{12} E^{2'}
\end{eqnarray*}
These pairs of sections are smooth on $U$ and independent on $U \setminus H_0$.   We compute from (\ref{eqn.E.action}) and (\ref{eqn.Estar.action})
\begin{eqnarray*}
z^t_{\mathcal{E}} (v(X)) = (1+tx_{11}) \cdot v(z^t.X) &  & z^t_{\mathcal{E}}(\iota(X)) = (1+tx_{11}) \cdot \iota(z^t.X) \\
(z^{-t}_{\mathcal{E}})^* (\tilde{v}(X)) = \tilde{v}(z^t.X) &  & (z^{-t}_{\mathcal{E}})^* (\tilde{\iota}(X)) = \tilde{\iota}(z^t.X) 
\end{eqnarray*}

Now define sections of $\mathcal{F}$ and $\mathcal{F^*}$ for $i = 2, \ldots, n$ by
\begin{eqnarray*}
w(X) = x_{11} E_1 + \cdots + x_{n1} E_n &  & \kappa_i(X) = E_i \\
\tilde{w}(X) = E^1 & &  \tilde{\kappa}^i(X) = x_{11}E^i - x_{i1}E^1
\end{eqnarray*}

These sections transform, for $i = 2, \ldots, n$, according to (\ref{eqn.F.action}) and (\ref{eqn.Fstar.action}) by
\begin{eqnarray*}
z^t_{\mathcal{F}}(w(X)) = w(z^t.X) & & z^t_{\mathcal{F}}(\kappa_i(X)) = \kappa_i(z^t.X) \\
(z^{-t}_{\mathcal{F}})^*(\tilde{w}(X)) = (1+tx_{11}) \cdot \tilde{w}(z^t.X) & & (z^{-t}_{\mathcal{F}})^*(\tilde{\kappa}^i(X)) = (1+tx_{11})\cdot \tilde{\kappa}^i(z^t.X)
\end{eqnarray*}

\subsection{Invariant section of the endomorphism bundle}
\label{sec.endo.section}

Now consider the sections
$$ \varphi' = v \otimes \tilde{\iota} \qquad \mbox{and} \qquad \varphi_i = \tilde{\kappa}^i \otimes w$$
of $\End \mathcal{E}^*$ and $\End \mathcal{F}$, respectively, for $i=2, \ldots, n$.  These preserve the decompositions (\ref{eqn.decomps.EF}).
They are each nilpotent endomorphisms of order two for any $X$: $(\varphi'_X)^2 = 0$ and $((\varphi_i)_X)^2 = 0$ for all $i$; in particular, they are trace-free.


The tensor $\varphi' \otimes \varphi_i$ is a section of the subbundle $\mbox{End}_0 \ \mathcal{E}^* \otimes \mbox{End}_0 \ \mathcal{F} \subset \End TU$, as in Section \ref{sec.inf.def}, corresponding to nontrivial deformations of the structure.  The flow acts on this section by
$$ (z^t)_*(\varphi'(X) \otimes \varphi_i(X)) = (1+tx_{11})^2 \cdot \varphi'(z^t.X) \otimes \varphi_i(z^t.X)$$

Define $q(X) = x_{12}^2 + x_{11}^2 + \cdots + x_{n1}^2$.  Note that $q(z^t.X) = (1+tx_{11})^{-2} \cdot q(X)$.  Then the section 
$$\Phi_i = \frac{1}{q} \varphi' \otimes \varphi_i$$
is $z^t$-invariant.  The coefficients of the components of $\Phi$ are rational functions in the variables $x_{12}, x_{11}, \ldots, x_{n1}$ with numerator degree four and denominator degree two.  They are smooth on $U \backslash SF$ and $C^1$ on $SF$, in particular at the origin.

Denote by $E^{i'}_{j'}$ the elementary endomorphism of $\mathcal{E}^*$ sending $E^{j'}$ to $E_{i'}$, and by $E_j^i$ the elementary endomorphism of $\mathcal{F}$ sending $E_i$ to $E_j$.  The coefficients of $\Phi_i$ are given by 
\begin{eqnarray*}
\frac{1}{q(X)} \cdot \left( -x_{11}x_{12}E_{1'}^{1'} - x_{12}^2 E_{1'}^{2'} + x_{11}^2 E_{2'}^{1'} + x_{11}x_{12} E_{2'}^{2'} \right) \\
\otimes \left( \sum_{k=1}^n x_{11}x_{k1} E_k^i - \sum_{k=1}^n x_{i1}x_{k1} E_k^1 \right) 
\end{eqnarray*}

which expands further, denoting $E^{i'}_{j'} \otimes E_k^\ell$ by $E^{i' \ell}_{j'k}$, as
\begin{eqnarray*}
\sum_k \frac{x_{12}x_{11}x_{i1}x_{k1}}{q(X)} E_{1'k}^{1'1} + \sum_k \frac{-x_{12}x_{11}^2x_{k1}}{q(X)} E_{1'k}^{1'i} + \sum_k \frac{x_{12}^2x_{i1}x_{k1}}{q(X)} E_{1'k}^{2'1} \\
+ \sum_k \frac{-x_{12}^2x_{11}x_{k1}}{q(X)} E_{1'k}^{2'i} + \sum_k \frac{-x_{11}^2x_{i1}x_{k1}}{q(X)} E_{2'k}^{1'1} + \sum_k \frac{x_{11}^3 x_{k1}}{q(X)} E_{2'k}^{1'i} \\
+ \sum_k \frac{-x_{12}x_{11}x_{i1}x_{k1}}{q(X)} E_{2'k}^{2'1} + \sum_k \frac{x_{12}x_{11}^2x_{k1}}{q(X)} E_{2'k}^{2'i}
\end{eqnarray*}

Of course, for any constant ${\bf c}=(c_2, \ldots, c_n)$, the endomorphism field
$$\Phi=\Phi_{\bf c} = \sum_{i=2}^n c_i \Phi_i$$
will be $z^t$-invariant and $C^1$.

\begin{proposition}
\label{main.prop}
 Let $\theta : TU \stackrel{\sim}{\rightarrow} \Hom(\BR^2,\BR^n)$ be the flat almost-Grass\-man\-nian structure on $U$.  Then for any ${\bf c} \neq {\bf 0}$, $(\Id + \Phi_{\bf c}) \circ \theta$ is a $\{ z^t \}$-invariant, $C^1$ almost-Grassmannian structure on $U$ not $C^1$ equivalent to $\theta$.
\end{proposition}

Fix ${\bf c} \neq {\bf 0}$ and denote $\Phi = \Phi_{\bf c}$.  We first show that $(\Id + \Phi) \circ \theta$ is an almost-Grassmannian structure on $U$.  Recall from above that $(\varphi')^2 = 0$ and $(\varphi_i)^2 = 0$ for all $i = 2, \ldots, n$; note that moreover, $\varphi_i \circ \varphi_j = 0$ for all $i,j = 2, \ldots, n$, so that $c_2 \varphi_2 + \cdots + c_n \varphi_n$ is also nilpotent of order two.  It follows that $\Phi_X$ is a nilpotent endomorphism of $(\mathcal{E}^* \otimes \mathcal{F})_X$ of order two for all $X \in U$, so the matrix exponential of $\Phi_X$ in $\SL(\mathcal{E}^* \otimes \mathcal{F})_X$ is simply $\Id + \Phi_X$.  We conclude that $\Id + \Phi_X$ is an isomorphism for all $X \in U$, so $(\Id + \Phi) \circ \theta$ is an almost-Grassmannian structure on $U$.  The $\{ z^t \}$-invariance holds by construction.

The derivatives of the rational coefficients in $\Phi$ are undefined on $SF$, the zero set of $q$.  The numerators are homogeneous polynomials in $x_{12}, x_{11}, \ldots, x_{n1}$ of degree five, with denominators all equal $q^2$.  Such functions extend continuously to $0$ on $SF$.

The final claim of the proposition is proved in the following section.

\subsection{Calculation of nonzero harmonic torsion terms}
\label{sec.proof.calc}

Recall from Section \ref{tors} that the harmonic torsion can be computed from any principal connection $\gamma \in \Omega^1(\Cal G_0,\lieg_0)$.  Such a connection is equivalent to a pair of (volume-compatible) linear connections $\nabla_{\Cal E^*}$ on $\Cal E^*$ and $\nabla_{\Cal F}$ on $\Cal F$. These induce a connection  on $\mathcal{E}^* \otimes \mathcal{F}$ of the form $\nabla_{\mathcal{E^*}} \otimes \Id_{\mathcal{F}} + \Id_{\mathcal{E^*}} \otimes \nabla_{\mathcal{F}}$. Via an almost-Grassmannian structure $\Psi : TU \stackrel{\sim}{\rightarrow} \mathcal{E^*} \otimes \mathcal{F}$ of class at least $C^1$, this connection can be pulled back to $TU$, and that pullback has a well-defined torsion. For the harmonic torsion,
we map back to $\Cal E^*\otimes\Cal F$ via $\Psi$ and project to the bundle associated to $ (S^2\BR^2\otimes\BR^{2*})_0\boxtimes (\La^2 \BR^{n*}\otimes\BR^n)_0,$ as in Section \ref{tors}.
Nonvanishing of the result is an obstruction to $C^1$ flatness.

 
As $\mathcal{E^*}$ and $\mathcal{F}$ are trivial bundles over $U$, we can use the trivial connections on each. By construction, the frame $\{E^{i'}_j\}$ of $\Cal E^*\otimes\Cal F$ is parallel.
Denote by $\nabla$ the pullback to $TU$ via our deformed almost-Grassmannian structure corresponding to $\Id + \Phi$.
The pullbacks of $\{E^{i'}_j\}$ comprise a framing of $TU$ by parallel vector fields $\{ \tilde E^{i'}_j \}$. The torsion is then determined by their brackets.
We compute a specific component of the torsion and apply Lemma \ref{harm-tors}
to show that the harmonic torsion is nonzero.

Consider the sections $E^{2'}_s$ for $s>1$ and $E^{2'}_1$ of $\Cal E^*\otimes\Cal F$. Both have rank one with kernel spanned by $E_{1'}$. Now $T(\tilde E^{2'}_s,\tilde E^{2'}_1)$ by construction depends only on the component of $T$ in the subbundle corresponding to the second summand in the decomposition \eqref{La2-decomp}. Consequently, Lemma \ref{harm-tors} implies nonvanishing of the harmonic torsion provided $T(\tilde E^{2'}_s,\tilde E^{2'}_1)(E_{1'})$ is not contained in the span of $E_1$ and $E_s$.
From the explicit description of $\Phi$ from above, we obtain for $s>1$
$$ \tilde{E}^{2'}_s = \theta^{-1} \circ (\Id + \Phi)^{-1}(E^{2'}_s) = \partial^{2'}_s - c_s \cdot \sum_{p'=1, k = 1}^{p'=2, k=n} \frac{x_{1p'} x_{11}^2 x_{k1}}{q(X)} \partial^{p'}_k,$$
while  
$$ \tilde{E}^{2'}_1 = \theta^{-1} \circ (\Id + \Phi)^{-1}(E^{2'}_1) = \partial^{2'}_1 + \sum_{p'=1,k=1,i=2}^{p'=2,k=n,i=n} c_i \cdot \frac{x_{1p'} x_{11} x_{k1} x_{i1}}{q(X)} \partial^{p'}_k.$$
Since these fields are parallel, the torsion is given by  
$$T(\tilde{E}^{2'}_s, \tilde{E}^{2'}_1) = -[\tilde{E}^{2'}_s, \tilde{E}^{2'}_1]=- \frac{c_s x_{11}^2}{q(X)} \cdot \sum_k \left( x_{k1} \partial^{2'}_k - \frac{2 x_{k1} x_{12}}{q(X)} \sum_{p'} x_{1p'} \partial^{p'}_k \right).$$

Mapping this vector field $\tilde{D}$ to $\Cal E^*\otimes \Cal F$ via $(\Id + \Phi)\circ\theta$, we obtain, in order of increasing \emph{net degree}---degree of numerator minus degree of denominator,
\begin{eqnarray*}
D := (\Id + \Phi) \circ \theta (\tilde{D})
= & - \frac{c_s x_{11}^2}{q(X)} \cdot \sum_k \left( x_{k1} E^{2'}_k - \frac{2 x_{k1} x_{12}}{q(X)} (x_{11} E^{1'}_k + x_{12} E^{2'}_k )  \right. \\
+ & \left.  x_{k1} \Phi(E^{2'}_k) - \frac{2 x_{k1} x_{12}}{q(X)} \left( (x_{11} \Phi(E^{1'}_k) + x_{12} \Phi(E^{2'}_k )  \right) \right)
\end{eqnarray*}
This is a continuous section of $\Cal E^*\otimes\Cal F$. Now compute
$$
D(E_{1'}) = \frac{c_s x_{11}^2}{q(X)} \cdot  \sum_k \frac{-2 x_{k1} x_{12} x_{11}}{q(X)} E_k + \ \mbox{higher-order terms},$$
where the higher-order terms have net degree at least three. This has nontrivial projection modulo $\mbox{span} \{E_1, E_s \}$ on an open, dense subset of any neighborhood of $0$ in $U$, provided $c_s \neq 0$.  

We conclude that the harmonic torsion of the deformed structure given by $\Id + \Phi$ is nontrivial on an open, dense subset of any neighborhood of 0.  This structure is thus inequivalent to $\Gr(2,n)$ on any open subset containing $0$ in its closure.

\subsection{Vanishing of harmonic curvature for $C^2$ structures}
\label{sec.C2.curvature}
Let $(\Cal G_0 \to M, \theta)$ be a $(2,n)$ almost-Grassmannian structure of regularity $C^2$ admitting a strongly essential flow $\{ z^t \}$ with higher-order fixed point $x_0$. 
We verify here that the proof of \cite[Thm 3.1]{cap.me.parabolictrans} applies to such a structure, so it has vanishing harmonic curvature on an open set containing the higher order fixed point in its closure.  

Let $(p: \Cal G \to M, \omega)$ be the $C^1$ prinicpal $P$-bundle and $C^0$ Cartan connection, respectively, given by Proposition \ref{prop.cg.kgeq2}.  There is a $C^1$ exponential map $\exp : \Cal G \times \lieg \rightarrow \Cal G$, giving $C^1$ exponential curves $\tilde{\gamma}_X(s)= \exp(u,sX)$ as in \cite[Def 1.4]{cap.me.parabolictrans}.   This differentiability is sufficient to apply all of the holonomy calculations of \cite[Sec 2]{cap.me.parabolictrans}.  

Let $\xi$ be the vector field generating $\{ z^t \}$.  For any $u \in p^{-1}(x_0)$, the value $\omega_{u}(\xi) \in \lieg_1$ (\cite[Sec 1.2]{cap.me.parabolictrans}); let $Z$ be the value for a particular choice of $u$.  The rank of $Z$ as an element of $\Hom(\BR^n, \BR^2)$ can be two or one.
 In either case, $Z$ defines a subset $\Cal T(X) \subset \lieg_{-1}$ comprising elements generating an $\mathfrak{sl}_2$-triple $\{X, A = [Z,X], Z \}$ (see \cite[Def 2.11]{cap.me.parabolictrans}).  Along exponential curves $\tilde{\gamma}_X$ for $X \in \Cal T(Z)$, the harmonic curvature and harmonic torsion must belong to the \emph{stable subspaces} $\mathbb{K}^{[st]}$ and $\mathbb{T}^{[st]}$, respectively, determined by $A$ (\cite[Def 2.13]{cap.me.parabolictrans}).  This restriction appears in Corollary 2.14 (1) of \cite{cap.me.parabolictrans}, which in turn rests on Proposition 2.9 of the same; the required property here is that the harmonic curvature and torsion are given by \emph{continuous} $P$-equivariant functions on $\Cal G$.

When $\rk Z = 2$, then $| \Cal T(Z) | = 1$, and $\mathbb{K}^{[st]} = 0$.  The harmonic curvature vanishes not only along the curve $\gamma_X = p \circ \tilde{\gamma}_X$, but on a neighborhood of $\gamma_X \backslash \{ x_0 \}$, as given by \cite[Prop 3.3]{cap.me.parabolictrans}.

When $\rk Z = 1$, as in the examples constructed above, then the \emph{strongly stable subspace} $\mathbb{K}^{[ss]}$ determined by $A$ is trivial.  Together with other purely algebraic features of the $\mathfrak{sl}_2$-triple and its representation on $\mathbb{K}$, this property suffices to again prove vanishing of the harmonic curvature on a neighborhood of $\gamma_X \backslash \{ x_0 \}$, as shown in \cite[Prop 3.5 (4a)]{cap.me.parabolictrans}.

The proof of \cite[Thm 1.3]{mn.1graded} that $(2,n)$-almost-Grassmannian manifolds with strongly essential automorphisms have vanishing harmonic torsion on an open set containing the higher-order fixed point in its closure is more involved, and requires higher regularity; in particular, the full Cartan curvature and a \emph{fundamental derivative} (see \cite[Sec 3.4]{mn.1graded}) of it must be continuous.  Thus it remains unclear whether a deformation with the properties in the conclusion of Theorem \ref{main.thm} could have higher regularity than $C^1$.

\subsection{Nonvanishing harmonic curvature on smooth locus}
\label{sec.nonzero.curv}

Our $C^1$ deformed structures $(\Id + \Phi) \circ \theta$ do not have a well-defined harmonic curvature tensor everywhere, but we can prove that the harmonic curvature for these structures is nontrivial on the smooth locus $U \backslash SF$, for sufficiently small values of ${\bf c}$. 
The \emph{infinitesimal change} of harmonic curvature produced by the infinitesimal deformation $\Phi_{\bf c}$ is the derivative at $t=0$ of the change of harmonic curvature produced by the deformations $(\Id + t \Phi_{\bf c}) \circ \theta = (\Id + \Phi_{t {\bf c}}) \circ \theta$.  For our aim, it suffices to show that the infinitesimal change of harmonic curvature caused by $\Phi_{\bf c}$ is nontrivial for ${\bf c} \neq {\bf 0}$. 

According to Theorem 3.6 of \cite{cap.deformations}, the infinitesimal change of harmonic curvature induced by an infinitesimal deformation can be computed on the smooth locus with the BGG sequence constructed from a certain linear connection on the adjoint tractor bundle. The operators in that BGG sequence act between sections of bundles associated to Lie algebra homologies $H_*(\lieg_1,\lieg)$, which are isomorphic to the Lie algebra cohomology spaces $H^*(\lieg_{-1},\lieg)$. These are representations of $\lieg_0$, which can be computed explicitly using Kostant's version of the Bott-Borel-Weil Theorem. The specific calculations for the Grassmannian case can be found in Section 3.5 of \cite{cap.soucek.subcomplexes} and 
in Section 4.1.3 Step (D) of \cite{cap.slovak.book.vol1}. 

In degree one, this representation is irreducible and isomorphic to $\mathfrak{sl}(\BR^{2*}) \boxtimes \mathfrak{sl}(\BR^n)$. 
In degree two, there are two irreducible components, one of which is the module $\mathbb T$ of Section \ref{tors} corresponding to the harmonic torsion. The other component $\mathbb{K} \subset \Lambda^2 \lieg_{-1}^* \otimes \lieg_0$ is the harmonic curvature module of Section \ref{sec.curv.and.connxn}.  It is the component of maximal highest weight in  
$$ (\Lambda^2 \BR^2 \boxtimes S^2 \BR^{n*}) \otimes (\mathfrak{sl}(\BR^{2*}) \oplus \mathfrak{sl}(\BR^n))$$
(recall the decomposition in (\ref{La2-decomp})), which turns out to be 
$$\mathbb{K} \cong \Lambda^2 \BR^2 \boxtimes (S^3 \BR^{n*} \otimes \BR^n)_0$$

The construction of BGG sequences provides an invariant differential operator $D$ mapping sections of 
$\mbox{End}_0 (\mathcal{E}^*) \boxtimes \mbox{End}_0 (\mathcal{F})$, corresponding to infinitesimal deformations, to sections of the bundle associated to $\mathbb K$, yielding the infinitesimal change of harmonic curvature (see Section 3.6 of \cite{cap.deformations}). 
  Any BGG operator admits a universal formula in terms of any \emph{distinguished connection} of the structure, its curvature and torsion, and their covariant derivatives. 
  In our case, the initial structure is the flat structure on the open set $U \backslash SF \subset \mbox{Gr}(2,n)$, for which the flat connection $\nabla_0$ induced by the trivial connections on $\Cal E^*$ and $\Cal F$ as in Section \ref{sec.proof.calc} is a distinguished connection.

Representation theory implies that $D$ must be a second-order operator.
Since $\nabla_0$ is torsion-free and flat, the universal formula for $D$ can only consist of applying two covariant derivatives, which automatically are symmetric, followed by tensorial operations induced by $\lieg_0$-equivariant maps on the inducing representations. The latter is a $\lieg_0$-equivariant map
$$
\rho : S^2(\lieg_{-1}^*)\otimes (\mathfrak{sl}(\BR^{2*}) \boxtimes \mathfrak{sl}(\BR^n))\to\mathbb K. 
$$
Similarly as in \eqref{La2-decomp}, we can decompose
$$
S^2(\lieg_{-1}^*)\cong (S^2\BR^2\boxtimes S^2\BR^{n*})\oplus (\La^2\BR^2\boxtimes\La^2\BR^{n*}),
$$
It is easy to see from representation theory that $\rho$ factors through the first summand. 
There is a unique homomorphism $S^2\BR^2\otimes(\BR^2\otimes\BR^{2*})_0\to\La^2\BR^2$ up to scale: the unique nonzero contraction with values in $\BR^2 \otimes \BR^2$, followed by an alternation.
There is also a unique homomorphism
$$S^2\BR^{n*}\otimes(\BR^n\otimes\BR^{n*})_0\to(S^3 \BR^{n*} \otimes \BR^n)_0$$
 up to scale: symmetrization of the three $\BR^{n*}$ components followed by projection on the module of trace-free elements. 


Using similar index notation as before, we now form the second derivative $\{ \nabla_{i'}^j \nabla_{\ell'}^m \Phi^{p' o}_{q' r} \}$, and then projection to $\mathbb{K}$ is achieved by 
 \begin{enumerate}
 \item contracting the indices $p'$ and $\ell'$;
 \item skew-symmetrizing the indices $i'$ and $q'$;
 \item symmetrizing the indices $j, m$, and $o$; and
 \item removing the trace of $r$ with $(j m o)$
 \end{enumerate}

Like in Section \ref{sec.proof.calc}, it now suffice to apply the operations in the $\BR^2$-part and then find a nonzero component which cannot be contained in the pure trace component. Namely, we now compute the term $\kappa_{2' 1' r}^{1 1 1}$, for $r > 1$, which evidently has this property.
\begin{eqnarray*}
 \kappa^{111}_{2' 1' r} & = & \frac{1}{2} \left(  \nabla^1_{2'} \nabla_{1'}^1 \Phi^{1' 1}_{1' r} + (\nabla_{2'}^1)^2 \Phi^{2' 1}_{1' r} - (\nabla^1_{1'})^2 \Phi^{1' 1}_{2' r} -  \nabla_{1'}^1 \nabla_{2'}^1 \Phi_{2' r}^{2' 1}  \right) \\
 & = & \frac{1}{2} \left( 2 \nabla^1_{2'} \nabla_{1'}^1 \Phi^{1' 1}_{1' r} + (\nabla_{2'}^1)^2 \Phi^{2' 1}_{1' r} - (\nabla^1_{1'})^2 \Phi^{1' 1}_{2' r} \right) 
\end{eqnarray*}
using that $\nabla^1_{2'} \nabla_{1'}^1 \Phi^{1' 1}_{1' r} = -  \nabla_{1'}^1 \nabla_{2'}^1 \Phi_{2' r}^{2' 1}$ by trace-freeness in the primed indices and flatness of $\nabla$. Recall from Section \ref{sec.endo.section}
$$ \Phi^{1' 1}_{1' r}  = \sum_{i > 1} c_i  \frac{x_{12} x_{11} x_{i1} x_{r1}}{q(X)}  \qquad  \Phi^{2' 1}_{1' r} = \sum_{i > 1} c_i  \frac{x_{12}^2 x_{i1} x_{r1}}{q(X)} \qquad \Phi^{1' 1}_{2' r} = - \sum_{i > 1} c_i  \frac{ x_{11}^2 x_{i1} x_{r1}}{q(X)}$$

Then compute, writing $q = q(X)$,
\begin{eqnarray*}
\nabla^1_{2'} \nabla_{1'}^1 \Phi^{1' 1}_{1' r} & = & \sum_{i > 1} c_i \left( \frac{x_{i1} x_{r1}}{q} -  \frac{2(x_{11}^2 + x_{12}^2) x_{i1} x_{r1}}{q^2} +  \frac{8 x_{12}^2 x_{11}^2 x_{i1} x_{r1}}{q^3} \right)  \\
(\nabla_{2'}^1)^2 \Phi^{2' 1}_{1' r} & = & \sum_{i > 1} c_i \left( \frac{2 x_{i1} x_{r1}}{q} - \frac{10 x_{12}^2 x_{i1} x_{r1}}{q^2} + \frac{8 x_{12}^4 x_{i1} x_{r1}}{q^3} \right) \\
(\nabla^1_{1'})^2 \Phi^{1' 1}_{2' r} & = & - \sum_{i > 1} c_i \left( \frac{2 x_{i1} x_{r1}}{q} - \frac{10 x_{11}^2 x_{i1} x_{r1}}{q^2} + \frac{8 x_{11}^4 x_{i1} x_{r1}}{q^3}  \right)  
\end{eqnarray*}

Finally,
$$ \kappa_{2' 1' r}^{111} = \sum_{i > 1} c_i \left( \frac{3x_{i1} x_{r1}}{q} -  \frac{7(x_{11}^2 + x_{12}^2) x_{i1} x_{r1}}{q^2} +  \frac{4 (x_{11}^2 + x_{12}^2)^2 x_{i1} x_{r1}}{q^3} \right) \neq 0.$$

\bibliographystyle{amsplain}
\bibliography{karinsrefs}

\end{document}